\documentclass[11pt]{article}
\usepackage{mathptmx}
\usepackage{amssymb}
 \usepackage{mathtools}
 \usepackage{pifont}
 \usepackage{multicol}
\usepackage{mathrsfs}
\usepackage{color}
\usepackage{multicol}
\usepackage{amsthm}
\usepackage{float}
\usepackage{tikz}
\usepackage{tikz-cd}
\usetikzlibrary{arrows}
\usetikzlibrary{shapes.arrows}   
\usetikzlibrary{positioning}
\usepackage{mathrsfs}
\usepackage{amssymb}
\usepackage{mathtools}

\definecolor{1}{rgb}{1,0.2,0.3}
\definecolor{2}{rgb}{0.1,0.3,0.5}
\definecolor{3}{rgb}{1,1,0}
\definecolor{4}{rgb}{255,255,255}
\newtheorem{theorem}{Theorem}[section]

\newtheorem{conjecture}[theorem]{Conjecture}
\newtheorem{question}[theorem]{Question}

\theoremstyle{definition}

 \newcommand{\R}{\mathbb{R}}

\theoremstyle{remark}

\usepackage[total={5.5in,8.5in},text={4.8in,7.65in},centering]{geometry}

\begin{document}
\pagenumbering{gobble} 

\tikzset
{
  x=.2in,
  y=.2in,
}


{\bf \Large \centering BRANKO GR\"UNBAUM IN MANY DIMENSIONS \\}

\bigskip

{ \centering \Large \bf Matthew Kahle\\}
{ \centering \Large Ohio State University\\}
{ \centering \Large mkahle@math.osu.edu\\}


\section{Recollections}

Even though he was nearly 90 years old, I was still surprised and sad to hear that Branko Gr\"unbaum had passed away. I took courses from Branko as a graduate student at the University of Washington. I still have my $3$-ring binders with the typed notes he passed out in his courses on ``Polyhedral Geometry'' and ``Configurations of Points and Lines''. These notes are full of meticulous illustrations, as was his style. I also read with him for several weeks of independent study, studying ``Venn Diagrams'', and I especially enjoyed one-on-one mathematical conversations with him.

His office was overflowing with mathematical art, mostly polyhedral sculptures he had made out of colored cardboard and other materials. He brought different models to class nearly every day, and he would pass them around for us to play with while he lectured.

Branko had exquisite geometric taste and intuition. Like many great mathematicians and artists, he was also idiosyncratic. He was particularly fond of various notions of symmetry---as evidenced by his writings on symmetric configurations of points and lines, symmetric Venn diagrams, non-convex regular polyhedra, etc.
%

The humble offering of this article is just to collect a few of my favorite open problems of Branko in one place. I can not even guarantee that all of these problems are strictly or originally due to Branko, but I think they are all questions that he was interested in at one time or another. I hope to give an interesting sampling of his mathematical interests, even if the sample is extremely small compared to his prolific output, and even if it is biased by my personal tastes.

\section{In dimension $2$}

\subsection{Symmetric Venn diagrams}

Let $F= \{ C_1, C_2, \dots, C_n \}$ be a family of $n$ simple closed curves in the plane. 
We say that $F$ is a \emph{Venn diagram} if
each of the $2^n$ subsets $X_1 \cap X_2 \cap \dots \cap X_n$ is nonempty and connected. Here $X_i$ denotes either the interior or exterior of the curve $C_i$. See Ruskey \cite{Ruskey97} for a 1997 survey.

Venn himself proved that Venn diagrams exist for every $n$. A Venn diagram is said to be \emph{symmetric} if it is invariant under a rotation of $2 \pi / n$. Branko wondered about the existence of symmetric Venn diagrams.  A theorem of Henderson \cite{Henderson63, WW08} is that $n$ must be prime. A natural question is whether this necessary condition is also sufficient.

This was an open question for a few decades, until it was finally resolved by Griggs, Killian, and Savage in \cite{GKS04}. Their construction involves some beautiful combinatorics around ``chain decompositions'' of the Boolean lattice. 

There is still more to do in this area, however. A Venn diagram is said to be \emph{simple} if no more than two curves intersect at any point. 

\begin{question}
Do simple, symmetric Venn diagrams with $n$ curves exist for every prime number $n$?
\end{question}

See Ruskey, Savage, and Wagon for a 2006 survey \cite{RSW06}. More recently, Mamakani and Ruskey gave the first examples of simple, symmetric Venn diagrams with $n=11$ and $n=13$ curves \cite{MR14}. 

\subsection{Gr\"unbaum's edge-coloring conjecture}

The Four Color Theorem for planar graphs is equivalent to the statement that the dual graph of every triangulated $2$-sphere is $3$-edge colorable. This is essentially Tait's 1880 observation that the Four Color Theorem is equivalent to showing that every cubic bridgeless planar graph is $3$-edge colorable \cite{Tait1880}.

In 1968, Gr\"unbaum conjectured a beautiful generalization.

\begin{conjecture}
If $ G $ is a simple loopless triangulation of an orientable surface $S$, then the dual graph of $ G $ is 3-edge-colorable.
\end{conjecture}

I learned as I was preparing this article that Gr\"unbaum's edge-coloring conjecture unfortunately does not hold in general. Kochol gave counterexamples: for every $g \ge 5$ there is a triangulation of the genus $g$ surface whose dual graph is not $3$-edge colorable \cite{Kochol09}.

Still, this leaves the question open for genus $g$ with $1 \le g \le 4$. In particular, the following seems to be  open.

\begin{conjecture}
If $ G $ is a simple loopless triangulation of a torus, then the dual graph of $ G $ is 3-edge-colorable.
\end{conjecture}

\section{In dimension $3$}

Branko spoke reverently about Steinitz's Theorem, and I think he felt that it deserved to be better known and appreciated by modern mathematicians.\\

\noindent {\bf Steinitz's Theorem.} A simple graph $G$ (i.e.\ with no loops or multiple edges) is the $1$-skeleton of a $3$-polytope if and only if $G$ is planar and $3$-connected. \\

The $1$-skeleton of a triangulated $2$-dimensional sphere is an edge-maximal planar graph, and hence is $3$-connected. So Steinitz's theorem implies that a triangulated $2$-sphere can be embedded in $\mathbb{R}^3$ as the boundary of a convex polytope. In particular, every triangulated $2$-sphere is polyhedral in $\mathbb{R}^3$.

A $2$-dimensional simplicial complex is said to be \emph{polyhedral} if it admits an embedding in $\mathbb{R}^3$ with every vertex corresponding to a point, every edge corresponding to a straight line segment, and every triangle corresponding to a flat triangle contained in an affine plane. (It must also actually be a topological embedding, i.e.\ a continuous, injective, map on the geometric realization of the complex.)

\begin{question} \label{quest:torus} Is every triangulation of the $2$-dimensional torus polyhedral? 
\end{question}

There are some reasons to think that the answer is no. Brehm gave an example of a triangulated Mobius strip which is not polyhedral \cite{Brehm83}. More recently, Leopold gave examples of triangulations of the non-orientable surface of genus 5 which he proved do not even admit polyhedral immersions in $\mathbb{R}^3$ \cite{Leopold17}. 

Update: G\"unter Ziegler pointed out to me that this problem was solved in \cite{ABE-M07}. Archdeacon, Bonnington, and Ellis-Monaghan showed that the answer to Question \ref{quest:torus} is affirmative.

\section{In dimension $4$}

The Ham Sandwich Theorem says that any $d$ objects in $\mathbb{R}^d$ can be simultaneously bisected by a single hyperplane \footnote{We will omit technical definitions and be deliberately vague about  what constitutes an ``object''. But any ``nice enough'' measure should suffice. Finite Borel measures, or finite points sets, are examples that are typically considered. For a more careful treatment, see Matou\v{s}ek's book \cite{Matousek03}.}.

A related dissection theorem in the plane says that any object in $\R^2$ can be partitioned into $4$ parts of equal area with two lines---this theorem is a nice exercise, and we invite the reader to work it out for themselves.

It is also true that any object in $\mathbb{R}^3$ can be partitioned into $8$ equal parts with $3$ planes. This theorem is considerably trickier than the corresponding $2$-dimensional exercise. A proof can be found in Chapter 4 of Edelsbrunner's book \cite{Edelsbrunner87}.

Branko asked if any object in $\mathbb{R}^{d}$ can be partitioned into $2^d$ equal parts with $d$ hyperplanes \cite{Grunbaum60}.

It is tempting to believe it, but unfortunately the answer is no for $d \ge 5$. The following counterexample for $d=5$ is due to Avis \cite{Avis84}, and similar counterexamples work in higher dimensions.

Consider $32$ distinct points along the moment curve in $\R^5$, which is parameterized by $$t \to (t, t^2, t^3, t^4, t^5) \mbox{ for } t \in \mathbb{R}.$$
Points on the moment curve are in general position. So any hyperplane can only intersect the moment curve in at most $5$ points. Then five hyperplanes can only intersect the moment curve in a total of at most $25$ points. So five hyperplanes can only divide these $32$ points into at most $26$ subsets and can not separate them all! 

Since the mass-partition conjecture is true for $d \le 3$ and false for $d \ge 5$, this leaves only the case $d=4$.

\begin{question} \label{4d} Can every object in $\R^4$ be partitioned into $16$ equal parts with $4$ hyperplanes? In particular, can every set of $16n$ points in general position in $\mathbb{R}^4$ be separated into $16$ subsets of $n$ points, using 4 hyperplanes?
\end{question}

According to Florian Frick, there is reason to think that the answer is yes. Let $\Delta(j,k)$ denote the minimum dimension $d$ such that every set of $j$ objects in $\R^d$ can be simultaneously partitioned into $2^k$ equal pieces by $k$ hyperplanes.

\begin{conjecture}[Ramos's conjecture]
For every $j, k \ge 1$, $$\Delta(j,k) = \left\lceil \frac{2^k - 1}{k} j \right\rceil.$$
\end{conjecture}

The Ham Sandwich Theorem is the case $k=1$ and Gr\"unbaum's mass partition conjecture is the case $j=1$. Ramos's conjecture would imply that $\Delta(1,4) = 4$, which would imply that the answer to Question \ref{4d} is yes. See Blagojevi\'{c} et. al \cite{BPFHZ18} for some recent progress. As one special case in the paper, the authors prove that $\Delta(1,4) \le 5$. That is, any sufficiently nice object in $\R^5$ can be partitioned into $16$ equal measure pieces by $4$ hyperplanes. 

\section{In $4$ dimensions and higher}

A well-known inequality for simple planar graphs with $v$ vertices and $e$ edges is that
$$e \le 3v.$$
Branko asked in 1970 \cite{Grunbaum70} if there are natural generalizations of this inequality for $d$-dimensional simplicial complexes embeddable in $R^{2d}$. For a simplicial complex $S$, let $f_i(S)$ denote the number of $i$-dimensional faces. One possible generalization of $e \le 3v$ is the following conjecture, which may be due to Gr\"unbaum.

\begin{conjecture} [Gr\"unbaum--Kalai--Sarkaria] \label{conj:embed} For every $d \ge 1$ there is a constant $C_d$ such that $$f_d (S)\le C_d f_{d-1}(S)$$ for every $d$-complex $S$ which is embeddable in $\R^{2d}$.
\end{conjecture} 

Various formulations of this conjecture have apparently been discussed for some time---see, for example, Dey \cite{Dey93}. Work of Kalai and Sarkaria \cite{Kalai91,Sarkaria92} suggests a precise formulation, which would in particular give the best possible constant $C_d$ in every dimension. Gr\"unbaum wrote that this question is ``still'' open in his 1970 article, but so far I have not been able to find any earlier reference to it in the literature.

Another possible generalization of the inequality for planar graphs $ e \le 3v$ would be to show that

\begin{conjecture} \label{conj:verts} For every $d \ge 1$ there is a constant $C_d$ such that
$f_d \le C_d f_0^d$ 
\end{conjecture}

Until recently, this was open even for $2$-complexes embeddable in $\mathbb{R}^4$, and the best known bound seemed to be $$f_d \le C_d f_0^{d+1-1/3^d},$$
which follows from some extremal hypergraph results of Erd\H{o}s \cite{Erdos64}. See for example the discussion by Wagner about ``forbidden subcomplexes'' in \cite{Wagner11}. Parsa improved this to $$f_d \le C_d f_0^{d+1-1/3^{d-1}}$$ in 2018 \cite{Parsa18}.

Karim Adiprasito has posted a preprint \cite{Adiprasito19} apparently proving Conjectures \ref{conj:embed} and \ref{conj:verts}. The preprint claims to prove much more, including the well-known ``$g$-conjecture'' characterizing the $f$-vectors of simplicial spheres. In \cite{AS19}, Adiprasito and Steinmeyer apparently give a somewhat simpler proof of the Gr\"unbaum--Kalai--Sarkaria conjecture.

\bigskip

{\bf Acknowledgements.} I am grateful to Karim Adiprasito, Florian Frick, Uli Wagner, and G\"unter Ziegler for helpful comments.

\bigskip

\bibliographystyle{plain}
\nocite{*}
\bibliography{Grefs}

\end{document}